\documentclass[12pt]{article}
\usepackage{amsmath,amssymb}
\usepackage{color}
\usepackage{graphicx}
\def\R{\mathbb{R}}

\def\arcsinh{\mathrm{arcsinh}}

\def\arccosh{\mathrm{arccosh}}
\def\tam{\mathrm{tam}}
\def\tamh{\mathrm{tamh}}
\def\arctam{\mathrm{arctam}}
\def\arctamh{\mathrm{arctamh}}
\newtheorem{Theorem}{\hspace*{\parindent}Theorem}
\newtheorem{Lemma}{\hspace*{\parindent}Lemma}
\newtheorem{Corollary}{\hspace*{\parindent}Corollary}

\title{A look at generalized trigonometric functions as functions of their two parameters and further new properties}
\author{D.B.\:Karp$^{\rm a,b}$\footnote{Corresponding author. E-mail: D. Karp --\emph{dmitrika@hit.ac.il},
E.\:Prilepkina -- \emph{pril-elena@yandex.ru}}~~and
E.G.\:Prilepkina$^{\rm c}$
\\[10pt]
\small{\textit{$\phantom{1}^a$Department of Mathematics, Holon Institute of Technology}}
\\
\small{\textit{52 Golomb Street,  P.O. Box 305, Holon 5810201, Israel}}
\\[7pt]
\small{\textit{$\phantom{1}^b$School of Economics and Management, Far Eastern Federal University,}}
\\
\small{\textit{10 Ajax Bay, Russky Island, Vladivostok, 690922, Russia}}
\\[7pt]
\small{\textit{$\phantom{1}^c$Institute of Applied Mathematics,}}
FEBRAS,
\\
\small{\textit{7 Radio Street, Vladivostok,  690041, Russia}}}
\date{}

\def\R{\mathbb{R}}

\def\arcsinh{\mathrm{arcsinh}}

\def\arccosh{\mathrm{arccosh}}
\def\tam{\mathrm{tam}}
\def\tamh{\mathrm{tamh}}
\def\arctam{\mathrm{arctam}}
\def\arctamh{\mathrm{arctamh}}
\begin{document}
\maketitle

\begin{abstract} 
Investigation of the generalized trigonometric and hyperbolic functions containing two parameters has been a very active research area over the last decade.  We believe, however, that their monotonicity and convexity properties with respect to parameters have not been thoroughly studied. In this paper, we make an attempt to fill this gap.  Our results are not complete; for some functions, we manage to establish (log)-convexity/concavity in parameters, while for others, we only managed the prove monotonicity, in which case we present necessary and sufficient conditions for convexity/concavity.  In the course of the investigation, we found two hypergeometric representations for the generalized cosine and hyperbolic cosine functions which appear to be new.  In the last section of the paper, we present four explicit integral evaluations of combinations of generalized trigonometric/hyperbolic functions in terms of hypergeometric functions.
\end{abstract}

\bigskip

\bigskip


\bigskip


\bigskip
\section{Introduction and Preliminaries}
The class of functions nowadays known as  the generalized trigonometric and hyperbolic functions can be traced back to the 
1879 paper of Lundberg \cite{Lundberg (1879)}. They were then rediscovered many times, typically in connection with the study of the “generalized circle” $|x|^p+|y|^p=1$ or as eigenfunctions of one-dimensional $p$-Laplacian.  Details can be found in the book \cite{Lang-book (2016)} by Lang and Edmunds. Earlier history has been  described by Peetre and Lindquist; see \cite{Lindqvist(1995), Lindqvist(2003)}.  The recently published textbook \cite{Poodiack(2022)}  makes the topic accessible to students.  The generalized trigonometric functions containing two parameters appeared in the works of \^{O}tani \cite{Otani(1984)} and Dr\'{a}bek and Man\'{a}sevich \cite{Drabek(1999)} in connection with $(p,q)$-Laplacian.  Namely, these authors defined a two-parameter generalized sine function $\sin_{p,q}(y)$ in \cite{Drabek(1999)} as the inverse function to the integral 
\begin{equation}\label{eq:arcsinpq}
	y=\arcsin_{p,q}(x)=\int_0^x\frac{dt}{(1-t^q)^{1/p}}=x{_{2}F_{1}}\left(\left.\!\!\begin{array}{c}1/p,1/q\\1+1/q\end{array}\right|x^q\!\right),
\end{equation}
with $p,q>1$, $x\in[0,1]$ and $y\in[0,\pi_{p,q}/2]$, where
$$
\frac{\pi_{p,q}}{2}=\int_0^1\frac{dt}{(1-t^q)^{1/p}}=\arcsin_{p,q}(1)=\frac{1}{q}B(1/q,1/p'),~~p'=p/(p-1).
$$
Here and below, ${}_2F_1$ stands for the Gauss hypergeometric function (\cite{Andrews(1999)}, Definition~2.1.5).
and $B(\cdot,\cdot)$ is Euler's beta function (\cite{Andrews(1999)}, Definition~1.1.3).  Note that $\pi_{p,q}\downarrow2$ as $q\uparrow\infty$ and as $p\uparrow\infty$ while monotonically decreasing from infinity.  
The  function $\arcsin_{p,q}:[0,1]\to[0,\pi_{p,q}/2]$ is an increasing homeomorphism and its inverse  
 $\sin_{p,q}(y)$ is well defined.  We can extend $\sin_{p,q}(y)$ to $[\pi_{p,q}/2,\pi_{p,q}]$ by $\sin_{p,q}(y)=\sin_{p,q}(\pi_{p,q}-y)$ and
 further to $[-\pi_{p,q},0]$ by oddness and to the whole $\mathbb{R}$ by $2\pi_{p,q}$ periodicity
 (\cite{Edmunds(2012)}, (2.4)). It is known that this extension is continuously differentiable on
$\R$ and $C^{\infty}$ everywhere except at the points $\{k \pi_{p,q}/2; k\in\mathbb{Z}\}$ (\cite{Edmunds(2012)}, p. 49).  The function $\cos_{p,q}(y)$ is then defined naturally by 
\begin{equation}\label{eq:cos-defined}
\cos_{p,q}(y)=\frac{d}{dy}\sin_{p,q}(y).
\end{equation}
It is not hard to show using \eqref{eq:arcsinpq} that this implies 
\begin{equation}\label{eq:cossinsum}
|\cos_{p,q}(y)|^p+|\sin_{p,q}(y)|^q=1    
\end{equation}
on the whole real line (\cite{Edmunds(2012)}, (2.7)) and we can drop the absolute values if we restrict our attention to the interval  $[0,\pi_{p,q}/2]$. Differentiating this identity with respect to $y$ and using  \eqref{eq:cos-defined} after a little rearrangement, we obtain:
\begin{equation}\label{eq:Dcos}
\frac{d}{dy}\cos_{p,q}(y)=-\frac{q}{p}\sin_{p,q}^{q-1}(y)\cos_{p,q}^{2-p}(y).
\end{equation}
It was observed in \cite{Baricz(2015)} that using the change of variable $u=(1-t^q)^{1/p}$ in
$$
\arccos_{p,q}(x)=\arcsin_{p,q}(w)=\int_0^w(1-t^q)^{-1/p}dt,\quad\text{where}\quad w=(1-x^p)^{1/q},
 $$
we obtain the representation
\begin{equation}\label{eq:arccospq}
	\arccos_{p,q}(x)=\frac{p}{q}\int_x^1(1-u^p)^{1/q-1}u^{p-2}du,\quad\text{where}\quad 0\leq x\le1.
\end{equation}

The generalized hyperbolic $(p,q)$-arcsine function is defined by 
\begin{multline}\label{eq:arcsinhpq}
	\arcsinh_{p,q}(x)\!=\!\!\int_0^x\frac{dt}{(1+t^q)^{1/p}}\!=\!x{_{2}F_{1}}\left(\left.\!\!\begin{array}{c}1/p,1/q\\1+1/q\end{array}\right|\!-\!x^q\!\right)
 \\
 =\bigg(\!\frac{x^p}{1+x^q}\!\bigg)^{\!\!1/p}\!\!\!
 {_{2}F_{1}}\left(\left.\!\!\begin{array}{c}1,1/p\\1+1/q\end{array}\right|\frac{x^q}{1+x^q}\!\right),
\end{multline}
where $p,q>1$, $x\in[0,\infty]$.  Denote
\begin{equation}\label{eq:hyperpi}
    \hat{\pi}_{p,q}=2\arcsinh_{p,q}(\infty)=
    \begin{cases}(2/q)B(1/p-1/q,1/q),& p<q
    \\
    \infty,& q\le{p}.
    \end{cases}
\end{equation}
Note that $\hat{\pi}_{p,q}\downarrow2$ as $q\uparrow\infty$ and $\hat{\pi}_{p,q}\downarrow\pi_{q}=2\pi/(q\sin(\pi/q))$ as $p\downarrow1$.  Miyakawa and Takeuchi observed in (\cite{Miyakawa(2021)}, p. 3) that $\hat{\pi}_{p,q}=\pi_{r,q}$ with $r=pq/(pq+p-q)$ and $r>1$ precisely when $p<q$.
The hyperbolic $(p,q)$-sine is then defined as the inverse function to $\arcsinh_{p,q}$ from \eqref{eq:arcsinhpq}:
$$
\sinh_{p,q}:[0,\hat{\pi}_{p,q}/2]\to[0,\infty].
$$
Naturally, 
\begin{multline}\label{eq:coskh-defined}
	\cosh_{p,q}(y)=\frac{d}{dy}\sinh_{p,q}(y)=\frac{1}{\frac{d}{dx}\int_0^x(1+t^q)^{-1/p}dt}_{\left|{x=\sinh_{p,q}(y)}\right.}
	\\
	=(1+x^q)^{1/p}_{\left|{x=\sinh_{p,q}(y)}\right.}=\left(1+[\sinh_{p,q}(y)]^q\right)^{1/p},
\end{multline}
so that 
$$
[\cosh_{p,q}(y)]^p-[\sinh_{p,q}(y)]^q=1.
$$
Differentiating this identity with respect to $y$ and using  \eqref{eq:coskh-defined} after a little rearrangement yields:
\begin{equation}\label{eq:Dcosh}
\frac{d}{dy}\cosh_{p,q}(y)=\frac{q}{p}\sinh_{p,q}^{q-1}(y)\cosh_{p,q}^{2-p}(y).
\end{equation}
If $x=\sinh_{p,q}(y)$ and $z=\cosh_{p,q}(y)$, then $x=(z^p-1)^{1/q}$. Hence,  substituting $x=(z^p-1)^{1/q}$ in \eqref{eq:arcsinhpq}, we obtain the formula  
$$
\arccosh_{p,q}(z)=\arcsinh_{p,q}(x)=\int_0^x(1+t^q)^{-1/p}dt,\quad\text{where}\quad x=(z^p-1)^{1/q}.
$$
Now, making the change of variable 
$$
u=(1+t^q)^{1/p},~~t=(u^p-1)^{1/q},~~dt=(p/q)(u^p-1)^{1/q-1}u^{p-1}du
$$ 
with endpoint correspondence  $t=0\Leftrightarrow{u=1}$, $t=y\Leftrightarrow{u=z}$, we can rewrite it as (by changing back $z\to{x}$)
\begin{equation}\label{eq:arccoshpq}
	\arccosh_{p,q}(x)=\frac{p}{q}\int_1^x(u^p-1)^{1/q-1}u^{p-2}du,\quad\text{where}\quad 1\le x<\infty.
\end{equation}
Note that $\arccosh_{p,q}(\infty)=\hat{\pi}_{p,q}/2$ as defined in \eqref{eq:hyperpi}, so that $\cosh_{p,q}(y)$ is defined on $[0,\hat{\pi}_{p,q}/2]$ with values in $[1,\infty]$.

The following hypergeometric representations seem to be new.
\begin{Lemma}\label{lm:arccos2F1}
		The following identities hold:
		\begin{equation}\label{eq:arccos2F1}
			\arccos_{p,q}(x)=x^{p-1}(1-x^p)^{1/q}
			{_{2}F_{1}}\left(\left.\!\!\begin{array}{c}1,1+1/q-1/p\\1+1/q\end{array}\right|1-x^p\!\right),~~0\le{x}\le1,
		\end{equation}
		\begin{equation}\label{eq:arccosh2F1}
			\arccosh_{p,q}(x)=\frac{1}{x}(x^p-1)^{1/q}
			{_{2}F_{1}}\left(\left.\!\!\begin{array}{c}1,1/p\\1+1/q\end{array}\right|1-x^{-p}\!\right),~~1\le{x}<\infty.
		\end{equation}
	\end{Lemma}

\textbf{Proof.} Indeed, we obtain from \eqref{eq:arccospq} using substitution $t=u^p$:
	\begin{multline*}
		\frac{q}{p}\arccos_{p,q}(x)=\frac{1}{p}\int_0^1(1-t)^{1/q-1}t^{-1/p}dt-\frac{1}{p}\int_0^{x^p}(1-t)^{1/q-1}t^{-1/p}dt
		\\
		=\frac{1}{p}B(1/q,1-1/p)-\frac{x^{p-1}}{p-1}
		{_{2}F_{1}}\left(\left.\!\!\begin{array}{c}1-1/p,1-1/q\\2-1/p\end{array}\right|x^p\!\right),
	\end{multline*}
where we used the easily verifiable formula (derived by changing the variable in Euler's representation (\cite{Andrews(1999)}, Theorem~2.2.1) by direct term-wise integration of the binomial expansion) 
\begin{equation}\label{eq:betaincomp}
	\int_{0}^{x}t^{a-1}(1-t)^{b-1}dt=\frac{x^a}{a}{_{2}F_{1}}\left(\left.\!\!\begin{array}{c}a,1-b\\a+1\end{array}\right|x\!\right).
\end{equation}
Next, we apply the connection formula (\cite{Andrews(1999)}, 2.3.13), which for the values of parameters in question reduces to (after application of the binomial theorem ${}_1F_0(a;-;x)=(1-x)^{-a}$)
	\begin{multline*}
		{_{2}F_{1}}\left(\left.\!\!\begin{array}{c}1-1/p,1-1/q\\2-1/p\end{array}\right|x^p\!\right)
		\\
		=x^{1-p}\frac{\Gamma(2-1/p)\Gamma(1/q)}{\Gamma(1+1/q-1/p)}-q\Big(1-\frac{1}{p}\Big)(1-x^p)^{1/q}{_{2}F_{1}}\left(\left.\!\!\begin{array}{c}1,1+1/q-1/p\\1+1/q\end{array}\right|1-x^p\!\right).
	\end{multline*}
Substituting this into the above expression and simplifying, we arrive at 
	\eqref{eq:arccos2F1}.
	
	To prove  \eqref{eq:arccosh2F1}, start with \eqref{eq:arccoshpq} and compute by changing variable $t=1/u^p$ and using \eqref{eq:betaincomp}:
		\begin{multline*}
		\frac{q}{p}\arccosh_{p,q}(x)=\frac{1}{p}\int_{x^{-p}}^{1}(1-t)^{1/q-1}t^{1/p-1/q-1}dt
		\\
		=
		\frac{1}{p}\int_{0}^{1}(1-t)^{1/q-1}t^{1/p-1/q-1}dt-\frac{1}{p}\int_{0}^{x^{-p}}(1-t)^{1/q-1}t^{1/p-1/q-1}dt
  \\
  =\frac{1}{p}B(1/q,1/p-1/q)-\frac{x^{-p(1/p-1/q)}}{p(1/p-1/q)}{_{2}F_{1}}\left(\left.\!\!\begin{array}{c}1/p-1/q,1-1/q\\1/p-1/q+1\end{array}\right|x^{-p}\!\right).
	\end{multline*}
Connection formula (\cite{Andrews(1999)}, 2.3.13) gives (again in view of ${}_1F_0(a;-;x)=(1-x)^{-a}$):
		\begin{multline*}
		{_{2}F_{1}}\left(\left.\!\!\begin{array}{c}1/p-1/q,1-1/q\\1/p-1/q+1\end{array}\right|x^{-p}\!\right)
		\\
		=\frac{\Gamma(1/p-1/q+1)\Gamma(1/q)}{\Gamma(1/p)}x^{1-p/q}-q(1/p-1/q)(1-x^{-p})^{1/q}
			{_{2}F_{1}}\left(\left.\!\!\begin{array}{c}1,1/p\\1+1/q\end{array}\right|1-x^{-p}\!\right).
	\end{multline*}
Substituting this into the above expression yields \eqref{eq:arccosh2F1}. $\hfill\square$

Takeuchi \cite{Takeuchi(2016),Miyakawa(2021),Miyakawa(2022)}  defined the following two-parameter generalizations of the tangent and hyperbolic tangent functions
\begin{equation}\label{eq:tam-defined}
\tam_{p,q}(y)=\frac{\sin_{p,q}(y)}{[\cos_{p,q}(y)]^{p/q}},~~~\tamh_{p,q}(y)=\frac{\sinh_{p,q}(y)}{[\cosh_{p,q}(y)]^{p/q}}
\end{equation}
and observed that their various properties are closer to those of the classical tangent and hyperbolic tangent functions than the properties of the functions 
$$
\tan_{p,q}(y)=\frac{\sin_{p,q}(y)}{\cos_{p,q}(y)},~~~ \tanh_{p,q}(y)=\frac{\sinh_{p,q}(y)}{\cosh_{p,q}(y)}
$$ 
defined in \cite{Edmunds(2012)}.  Note also that $\tam_{p,p}(y)=\tan_{p,p}(y)$ and $\tamh_{p,p}(y)=\tanh_{p,p}(y)$. We further contribute to the justification of this viewpoint by observing that the inverse functions $\arctam_{p,q}$, $\arctamh_{p,q}$ possess simple integral representations permitting application of Lemma~\ref{lm:derivatives} below to study monotonicity and convexity with respect to parameters. An integral representation of $\tam_{p,q}(y)$ is obtained as follows. Dividing relation \eqref{eq:cossinsum} restricted to $[0,\pi_{p,q}/2]$ by $\cos_{p,q}^p(y)$, we can cast it to the form:
\begin{equation}\label{eq:tamcos}
1+\tam_{p,q}^{q}(y)=\cos_{p,q}^{-p}(y).
\end{equation}
Differentiation then yields:
\begin{multline*}
\frac{d}{dy}\tam_{p,q}(y)=\frac{\cos_{p,q}^{1+p/q}(y)+\sin_{p,q}^{q}(y)\cos_{p,q}^{1+p/q-p}(y)}{\cos_{p,q}^{2p/q}(y)}
\\
=\cos_{p,q}^{1-p/q}(y)(1+\tam_{p,q}^{q}(y))=(1+\tam_{p,q}^{q}(y))^{1/q-1/p+1},
\end{multline*}
where we applied  \eqref{eq:tamcos} in the last equality (alternatively, one can differentiate \eqref{eq:tamcos} directly). This differentiation formula was first found in (\cite{Miyakawa(2021)}, p. 16). Hence, for the inverse function, we obtain the representation
\begin{equation}\label{eq:arctamintegral}
\arctam_{p,q}(x)=\int_{0}^{x}\frac{dt}{(1+t^q)^{1+1/q-1/p}}.
\end{equation}
This representation is also seen from the connection formula $\tam_{p,q}(y)=\sinh_{r,q}(y)$, $r=pq/(pq+p-q)$ (\cite{Miyakawa(2022)}, Theorem~2.1).
In particular, writing $p'$ for the conjugate exponent found from the equation $1/p'+1/p=1$, we will have \cite{Takeuchi(2016)} [p.1006]:
\begin{equation}
\arctam_{p',p}(x)=\int_{0}^{x}\frac{dt}{(1+t^p)^{2/p}}.
\end{equation}

In a similar fashion, (\cite{Miyakawa(2022)}, Theorem~2.2) indicates that 
\begin{equation}
\tamh_{p,q}(y)=\sin_{r,q}(y),~~\text{where}~r=pq/(pq+p-q),
\end{equation}
so that from \eqref{eq:arcsinpq}, we have 
\begin{equation}\label{eq:arctamhintegral}
\arctamh_{p,q}(x)=\int_{0}^{x}\frac{dt}{(1-t^q)^{1+1/q-1/p}}
\end{equation}
and 
\begin{equation}\label{eq:tamhp}
\arctamh_{p',p}(x)=\int_{0}^{x}\frac{dt}{(1-t^p)^{2/p}}.
\end{equation}

The literature on the two-parameter trigonometric and hyperbolic functions has grown rapidly over the last two decades and we will not attempt to survey it in this introduction. The reader is invited to consult the paper \cite{Lindqvist(2004)} for a historical overview, the introductory parts of \cite{Karp(2015),Kobayashi(2019),Takeuchi(2019)}, the survey \cite{Yin(2019)} and the books \cite{Lang-book (2016),Poodiack(2022)}  for numerous further references to more modern research and connections with various areas of mathematics. Inequalities for the generalized trigonometric functions with two parameters have been studied in \cite{Baricz(2015),Bhayo(2012)} for fixed parameters.  Research on convexity properties with respect to parameters is scarce. A few exceptions include \cite{BBV(2015)}, where Tur\'{a}n-type inequalities for the inverse trigonometric functions, also those with two parameters, were established, and our paper \cite{Karp(2015)} concerned with functions with one parameter. Both are listed in the survey paper  (\cite{Yin(2019)}, Section~5). The purpose of this note is two-fold. First, we investigate the monotonicity and convexity properties of the generalized trigonometric and hyperbolic functions with respect to their parameters $p$ and $q$.  We believe that these types of results have not previously appeared in the literature.  This is shown in the subsequent Section~\ref{sec2}. Second, we derive several integrals of the inverse and direct functions in terms of hypergeometric functions, which we also believe to be new. These results can be found in Section~\ref{sec3}. Finally, our Lemma~\ref{lm:arccos2F1} above gives new representations for the inverse generalized cosine and hyperbolic cosine functions, which are useful both from the theoretical viewpoint and for an effective computation of these functions and their inverses.

\section{Monotonicity and Convexity in Parameters}\label{sec2}

The following lemma will be our key tool for the forgoing investigation of the  properties of the generalized trigonometric functions. Its proof is an exercise in calculus; details can be found in (\cite{Karp(2015)}, Lemma~1).

\begin{Lemma}\label{lm:derivatives}
Suppose $I,J$ are convex subsets of $\R$.  Suppose $f(p,x)\in{C^2(J\times{I})}$ is strictly monotone on $I$ for each fixed $p\in{J}$ so that $y\to g(p,y):=f^{-1}(p,y)$ is well defined and monotone on $f(I)$ for
each fixed $p\in{J}$.  Then, the following relations hold true:
\begin{align}
&\frac{\partial}{\partial{p}}g(p,y)=-f'_{p}/f'_x,\label{eq:dgdalpha}
\\[5pt]
&\frac{\partial^2}{\partial{p^2}}g(p,y)=\left(f'_x\right)^{-2}
\left\{2f'_{p}f''_{xp}-f'_xf''_{pp}-\left(f'_{p}\right)^2f''_{xx}/f'_x\right\},
\label{eq:d2gd2alpha}
\\[5pt]
&\frac{\partial^2}{\partial{p^2}}\log[g(p,y)]=(xf'_x)^{-2}\left\{2xf'_{p}f''_{p{x}}-xf'_xf''_{pp}-x(f'_{p})^2f''_{xx}/f'_x-(f'_{p})^2\right\},
\label{eq:d2loggd2alpha}
\end{align}
where $x$ on the right is related to $y$ on the left by $y=f(p,x)$ or $x=g(p,y)$.
\end{Lemma}

{\bf{Remark}} 
Formulas \eqref{eq:d2gd2alpha} and \eqref{eq:d2loggd2alpha} can also be written in the following form:
$$
\frac{\partial^2}{\partial{p^2}}g(p,y)=\frac{1}{2}\frac{\partial}{\partial{x}}\left(\frac{f'_{p}}{f'_{x}}\right)^{\!\!2}-\frac{\partial}{\partial{p}}\left(\frac{f'_{p}}{f'_{x}}\right)
$$
and
$$
\frac{\partial^2}{\partial{p^2}}\log[g(p,y)]=\frac{1}{2x}\frac{\partial}{\partial{x}}\left(\frac{f'_{p}}{f'_{x}}\right)^{\!\!2}
-\frac{1}{x}\frac{\partial}{\partial{p}}\left(\frac{f'_{p}}{f'_{x}}\right)
-\left(\frac{f'_{p}}{xf'_x}\right)^2.
$$

We summarize the monotonicity properties in the following theorem. 
\begin{Theorem}\label{th:monotonicity}
The functions $p\to\sin_{p,q}(y)$ and $p\to\tam_{p,q}(y)$ are increasing on $(1,\infty)$, while the functions $p\to\sinh_{p,q}(y)$, $p\to\cosh_{p,q}(y)$ and $p\to\tamh_{p,q}(y)$ are decreasing in the same interval for each fixed $y\in (0,1)$ and $q\in(1,+\infty)$. The functions  $q\to\sin_{p,q}(y)$  and $q\to\tamh_{p,q}(y)$  are increasing on $(1,\infty)$ for each fixed $y\in (0,1)$ and $p\in(1,+\infty)$. The functions $p\to\cos_{p,q}(y)$, $q\to\cos_{p,q}(y)$, $q\to\cosh_{p,q}(y),$  $q\to\tam_{p,q}(y)$ and $q\to\sinh_{p,q}(y)$ are, generally speaking, not monotonic (i.e., for each of them, there exist values of $y$ such that no monotonicity takes place). 
\end{Theorem}

\textbf{Proof.}  Note that $\inf_{p,q>1}\pi_{p,q}=\inf_{p,q>1}\hat\pi_{p,q}=1$. Hence, for the values of the argument $y\in(0,1)$, all functions specified in the theorem are well defined on $p,q\in(1,\infty)$ and are given as the inverse functions to the integrals described in the introduction.

According to Lemma~\ref{lm:derivatives} and representation \eqref{eq:arcsinpq}, for $x<1$ we obtain 
$$
\frac{\partial}{\partial{p}}\sin_{p,q}(y)=-\frac{(1-x^q)^{1/p}}{p^2}\int_0^x(1-t^q)^{-1/p}\log(1-t^q)dt>0;
$$
$$
\frac{\partial}{\partial{q}}\sin_{p,q}(y)=-\frac{(1-x^q)^{1/p}}{p}\int_0^xt^q(1-t^q)^{-1-1/p}\log(t)dt>0.
$$

In a similar fashion, for $x>0$ from \eqref{eq:arcsinhpq} 
$$
\frac{\partial}{\partial{p}}\sinh_{p,q}(y)=-\frac{(1+x^q)^{1/p}}{p^2}\int_0^x(1+t^q)^{-1/p}\log(1+t^q)dt<0.
$$ 

We omit similar verification of the monotonicity of the functions $\tam_{p,q}(y)$ and $\tamh_{p,q}(y)$ with  the signs of derivatives obvious from the resulting expressions. Numerical tests show that the functions $p\to\cos_{p,q}(y)$, $q\to\cos_{p,q}(y)$, $q\to\cosh_{p,q}(z),$ $q\to\sinh_{p,q}(y)$, $q\to\cosh_{p,q}(y)$ are not monotonic. For example, Figure~1 shows that $q\to\cosh_{p,q}(y)$ is not a monotone function for $p=1.7$ and $y$ near $1$, as we see that $\cosh_{p,q_2}(y)<\cosh_{p,q_1}(y)$ when $y=0.88$, but $\cosh_{p,q_2}(y)>\cosh_{p,q_1}(y)$ when $y=0.99$. This can be made rigorous by computation with guaranteed precision at fixed values of the argument and the second parameter. For another example, the derivative
$$
\frac{\partial}{\partial{q}}\sinh_{p,q}(y)=-\frac{(1+x^q)^{1/p}}{p}\int_0^xt^q(1+t^q)^{-1-1/p}\log(t)dt
$$ 
changes sign with varying $q$ if we take $x=13/10$, $p=5$.

According to Lemma~\ref{lm:derivatives}, to prove that $p\to\cosh_{p,q}(x)$ is decreasing, we need to show that
$$
\frac{\partial}{\partial{p}}\cosh_{p,q}(y)=-f_{p}'/f_{x}'<0,
$$
where
$$
f(x)=\arccosh_{p,q}(x)=\int_0^{(x^p-1)^{1/q}}(1+t^q)^{-1/p}dt.
$$
Differentiation using the Leibniz integral rule yields ($x>1$):
$$
f_{p}'= \frac{1}{p^2}\int_0^{(x^p-1)^{1/q}}(1+t^q)^{-1/p}\log(1+t^q)dt+\frac{x^{p-1} ( x^p-1)^{1/q-1} \log(x)}{q}>0,
$$
$$ 
f_{x}'=\frac{p x^{ p-2} (x^p-1)^{1/q-1}}{q}>0. 
$$

\includegraphics[width=13.7 cm]{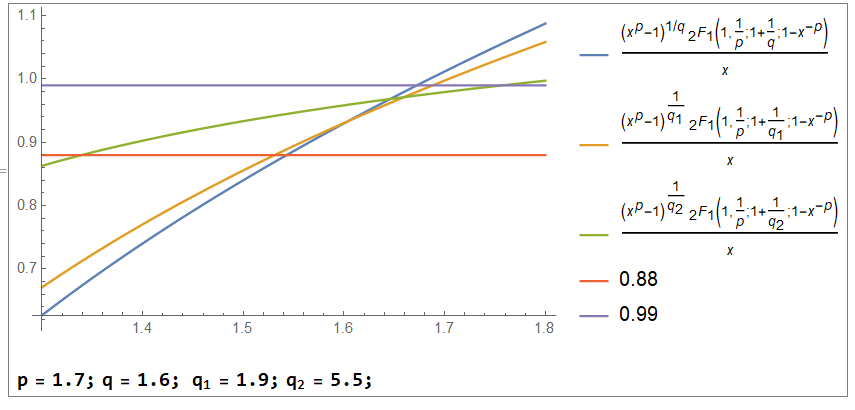}

{\small{\textbf{Figure~1.} The figure shows graphs of the hyperbolic $(p,q)$-arccosine at different values of $q$. By comparing the points of intersection of the graphs with the level lines $y=0.88$ and $y=0.99$, one can notice the absence of monotonicity  of the hyperbolic $(p,q)$-cosine as a function of $q$.}}  $\hfill\square$

\medskip

We proceed with some standard definitions.  A positive function $f$ defined on a finite or infinite interval $I$ is said to be logarithmically convex, or log-convex, if its logarithm is convex, or, equivalently,
$$
f(\lambda{x}+(1-\lambda)y)\leq[f(x)]^{\lambda}[f(y)]^{1-\lambda}~~\text{for all}~x,y\in{I}~~\text{and}~\lambda\in[0,1].
$$
The function $f$ is log-concave if the above inequality is reversed.  If the inequality is strict for $\lambda\in(0,1)$, then the respective property is called strict.  It can be seen from the definitions that log-convexity implies convexity and concavity implies log-concavity but not vice versa. 

In the context of generalized trigonometric and hyperbolic functions, the function $f(x,p)$ from Lemma~\ref{lm:derivatives} is given by an integral of a positive function $\varphi(t,p)$ of the form 
\begin{equation}\label{eq:gen1}
    y=f(x,p)=\int_{0}^{x}\varphi(t,p)dt
\end{equation}
which defines $x=g(y,p)$ implicitly for  $x\in I$.

\begin{Lemma}\label{lemma:conv}
  Suppose $x=g(p,y)$ is defined as the inverse to the integral of a  positive function $\varphi(t,p)$ of the form \eqref{eq:gen1}.  Then, convexity/concavity of  $p\to{g(y,p)}$ is determined by the sign of the~expression 
\begin{multline}\label{eq:conv1}
2\varphi'_p(x,p)\int_0^x\varphi'_p(t,p)dt-\varphi(x,p)\int_0^x\varphi''_{pp}(t,p)dt
\\
-\left(\int_0^x\varphi'_p(t,p)dt\right)^2\varphi'_x(x,p)/\varphi(x,p),
\end{multline}

while logarithmic convexity/concavity is determined by the sign of 
\begin{multline}\label{eq:logconv2}
2x\varphi'_p(x,p)\int_0^x\varphi'_p(t,p)dt-x\varphi(x,p)\int_0^x\varphi''_{pp}(t,p)dt-\\ (1+x\varphi'_x(x,p)/\varphi(x,p))\left(\int_0^x\varphi'_p(t,p)dt\right)^2.
\end{multline}

In particular, if for all $x\in I$ and $p>1$, the following inequalities hold:
 $$
 \varphi_x'(x,p)>0,\ D_1(x,p)=(\varphi'_p(x,p))^2-\varphi'_x(x,p)\int_0^x\varphi''_{pp}(t,p)dt<0,
 $$ 
 then the function $p\to g(p,y)$ is concave on $p>1$.  For logarithmic concavity, the following conditions are sufficient: 
$$
\varphi_x'(x,p)>0, D_2(x,p)=x(\varphi'_p(x,p))^2-\left(1+x\frac{\varphi'_x(x,p)}{\varphi(x,p)}\right)\varphi(x,p)\int_0^x\varphi''_{pp}(t,p)dt<0.
$$
\end{Lemma}
\textbf{Proof}. We have 
$$
f'_x=\varphi(x,p), \ f''_{xx}=\varphi'_x(x,p),\ f''_{xp}=\varphi'_p(x,p),\
$$
$$
f'_p=\int_0^x\varphi'_p(t,p)dt, \ f''_{pp}=\int_0^x\varphi''_{pp}(t,p)dt,
$$
and application of  Lemma~\ref{lm:derivatives} yields \eqref{eq:conv1} and \eqref{eq:logconv2}.

Next, denote $z=\int_0^x\varphi'_p(t,p)dt$. Verifying the concavity/logarithmic concavity, we deal with the quadratic forms 
\begin{equation}\label{eq:conv2}
-\frac{\varphi'_x(x,p)}{\varphi(x,p)} z^2+
2\varphi'_p(x,p) z-\varphi(x,p)\int_0^x\varphi''_{pp}(t,p)dt
\end{equation}
or
\begin{equation}\label{eq:logconv3}
-\left(1+x\frac{\varphi'_x(x,p)}{\varphi(x,p)}\right)z^2+
2x\varphi'_p(x,p) z-x\varphi(x,p)\int_0^x\varphi''_{pp}(t,p)dt.
\end{equation}
These quadratic forms are clearly negative for all $z$ if their respective discriminants and the leading coefficient are both negative.  $\hfill\square$

\begin{Corollary} \label{cr:psinpq}
The function $p\to\sin_{p,q}(y)$ is  concave \emph{(}strictly concave\emph{)} on $(1,\infty)$ for a fixed $y\in[0,1](=\cap_{p>1}[0,\pi_{p,q}/2])$ and
a fixed  $q\in(1,\infty)$ iff for all $p\in(1,\infty)$
\begin{multline}\label{eq:psinpq}
\frac{qx^{q-1}}{p(1-x^q)^{1-1/p}}\left(\int\limits_0^x\phi_1(p,t)dt\right)^2
-2\log(1-x^q)\int\limits_0^x\phi_1(p,t)dt
\\
+\int\limits_0^x\phi_{2}(p,t)dt\ge0\:(>0),
\end{multline}
where
$$
\phi_1(p,t)=\frac{\log(1-t^q)}{(1-t^q)^{1/p}},
~~~~\phi_2(p,t)=\frac{\log(1-t^q)(\log(1-t^q)-2p)}{(1-t^q)^{1/p}}.
$$
\end{Corollary}

\begin{Corollary} \label{cr:qsinpq}
The function $q\to\sin_{p,q}(y)$ is concave \emph{(}strictly concave\emph{)} on $(1,\infty)$ for a fixed  $y\in[0,1](=\cap_{q>1}[0,\pi_{p,q}/2])$ and a fixed
$p\in(1,\infty)$ iff for all $q\in(1,\infty)$
\begin{equation}\label{eq:qsinpq}
\frac{qx^{q-1}}{p(1-x^q)^{1-1/p}}\left(\int\limits_0^x\beta_1(q,t)dt\right)^2-\frac{2x^q\log x}{(1-x^q)}
\int\limits_0^x\beta_1(q,t)dt+\int\limits_0^x\beta_2(q,t)dt\ge0\:(>0),
\end{equation}
where
$$
\beta_1(q,t)=\frac{t^q\log(t)}{(1-t^q)^{1+1/p}},~~~\beta_2(q,t)=\frac{t^q(t^q+p)[\log(t)]^2}{(1-t^q)^{2+1/p}}.
$$
\end{Corollary}

Recall that $\pi_{q}=2\pi/(q\sin(\pi/q))$.
\begin{Corollary} \label{cr:psinhp}
The function $p\to\sinh_{p,q}(y)$ is log-convex \emph{(}strictly log-convex\emph{)}  on $(1,\infty)$ for a fixed $y\in[0,\pi_q](=\cap_{p>1}[0,\hat{\pi}_{p,q}])$ and
a fixed $q\in(1,\infty)$ iff for all $p\in(1,\infty)$
\begin{multline}\label{eq:psinhpq}
\frac{(q-p)x^q-p}{px(1+x^q)^{1-1/p}}\left(\int_0^x\lambda_1(p,t)dt\right)^2+2\log(1+x^q)\int_0^x\lambda_1(p,t)dt
\\
+\int_0^x\lambda_{2}(p,t)dt\ge0\:(>0),
\end{multline}
where
$$
\lambda_1(p,t)\!=\!\frac{\log(1+t^q)}{(1+t^q)^{1/p}},
~~~
\lambda_2(p,t)\!=\frac{(2p-\log(1+t^q))\log(1+t^q)}{(1+t^q)^{1/p}}.
$$
\end{Corollary}

\begin{Corollary} \label{cr:qsinhp}
The function $q\to\sinh_{p,q}(y)$ is log-convex \emph{(}strictly log-convex\emph{)} on $(1,\infty)$ for a fixed  $y\in[0,1](=\cap_{q>1}[0,\hat{\pi}_{p,q}])$ and a fixed $p\in(1,\infty)$ iff for all $q\in(1,\infty)$
\begin{multline}\label{eq:qsinhpq}
\frac{p+(p-q)x^q}{px(1+x^q)^{1-1/p}}\left(\int_0^x\alpha_1(q,t)dt\right)^2+\frac{2x^q\log x}{(1+x^q)}\int_0^x\alpha_1(q,t)dt
\\
+\int_0^x\alpha_{2}(q,t)dt\le0\:(<0),
\end{multline}
where
$$
\alpha_1(q,t)\!=-\frac{t^q\log(t)}{(1+t^q)^{1+1/p}},
	~~~
\alpha_{2}(q,t)\!=\frac{t^q(t^q-p)[\log(t)]^2}{(1+t^q)^{2+1/p}}.
$$
\end{Corollary}

The main results of this section are the following theorems.

\begin{Theorem}\label{th:tamhpp} 
The function $p\to\tamh_{p',p}(y)$, defined as the inverse to the integral  \eqref{eq:tamhp}, is logarithmically concave on half-line $p>1$. 
\end{Theorem} 

\textbf{Proof.}  According to  \eqref{eq:tamhp}, the kernel of the integral representation of the inverse function $\arctamh_{p',p}$ is given by $\varphi(x,p)=(1-x^p)^{-2/p}=:w(x,p)$. Differentiation then yields
$$
\varphi'_p(x,p)=2w(x,p)\left(\frac{x^p\mathrm{log}\, x}{p(1-x^p)}+\frac{\mathrm{log}\,(1-x^p)}{p^2}\right)=2w(x,p)\eta(x,p),
$$
where the second equality is the definition of $\eta(x,p)$. Further,
$$\varphi'_{x}(x,p)=2w(x,p)\frac{x^{p-1}}{1-x^p},\ 
\varphi''_{pp}(x,p)=\frac{(w'_p(x,p))^2}{w(x,p)}-2\frac{w'_p(x,p)}{p}+2w(x,p)\frac{x^p\log^2 x}{p(1-x^p)^2}, 
$$
$$
\varphi''_{pp}(x,p)=4w(x,p)\eta^2(x,p)-4w(x,p)\frac{\eta(x,p)}{p}+2w(x,p)\frac{x^p\log^2 x}{p(1-x^p)^2}.
$$

To establish logarithmic convexity  of $p\to\tam_{p',p}(x)$ by Lemma \ref{lemma:conv}, it is sufficient to check that  
$$ 
D_1(x)=4x\eta^2(x,p)\frac{1-x^p}{1+x^p}w(x,p)-\int_0^x\varphi''_{pp}(t,p)dt<0.
$$
Differentiating with respect to $x$ and substituting $u=x^p$, we obtain  
\begin{multline*}
D'_1(x)=-\frac{(1-u)^{-1-2/p}}{p^3(1+u)^2}
\Bigg\{2u(1-u)\log^2(u) 
-4(1-u)\log(1-u)
\\
\times\Big((1+u)^2-2u\log(1-u)\Big)-4u\log(u)\Big((1+u)^2+2(1-u)\log(1-u)\Big)\Bigg\}. 
\end{multline*}
Ascertaining by standard means that   
$(1+u)^2+2(1-u)\log(1-u)>0$ 
in the interval $u\in(0,1)$; we conclude that $D'_1(x)<0$, and as $D_1(0)=0$, this completes the proof of the logarithmic concavity of $p\to\tam_{p',p}(x)$. $\hfill\square$

\medskip

The following lemma is a guise of monotone L'H\^{o}spital rule \cite{Pinelis(2006)}.  See \cite{Karp(2015)} (Lemma~2) for a detailed proof.

\begin{Lemma}\label{lm:LHospital}
Suppose $f$, $g$ are continuously differentiable on a finite interval $(a,b)$,
$f(a)=g(a)=0$ and $gg'>0$ on $(a,b)$. If $f'/g'$ is decreasing on $(a,b)$, then $f/g>f'/g'$ on $(a,b)$.
\end{Lemma}

The above lemma will be applied in the proof of the following theorem.

\begin{Theorem}\label{th:sinpq}
For each fixed $y\in[0,1](=\cap_{p>1}[0,\pi_{p,q}/2])$ and $q>1$, the function
$p\to\sin_{p,q}(y)$ is concave on $(1,\infty)$.  For each fixed $y\in[0,1](=\cap_{q>1}[0,\pi_{p,q}/2])$ and $p>1$, the function $q\to\sin_{p,q}(y)$ is concave on $(1,\infty)$.
\end{Theorem}

\textbf{Proof.}  (i) According to Corollary \ref{cr:psinpq}, to prove concavity of $p\to\sin_{p,q}(x)$, we need to show (\ref{eq:psinpq}). Consider the quadratic
$$F(z)=\frac{qx^{q-1}}{p(1-x^q)^{1-1/p}}z^2
-2\log(1-x^q)z+\int\limits_0^x\frac{\log(1-t^q)(\log(1-t^q)-2p)}{(1-t^q)^{1/p}}dt.$$
We need to show that $F(z)>0$ for 
$$
z=\int_{0}^{x}\frac{\log(1-t^q)}{(1-t^q)^{1/p}}dt.
$$  
Since the coefficient at $z^2$ is positive, this quadratic is positive for all real $z$ if its discriminant $D$ is negative.  Compute
$$
D=\big(2\log(1-x^q)\big)^2-4\frac{qx^{q-1}}{p(1-x^q)^{1-1/p}}\int\limits_0^x\frac{\log(1-t^q)(\log(1-t^q)-2p)}{(1-t^q)^{1/p}}dt<0
$$
which is equivalent to
$$
G(x):=\frac{p(1-x^q)^{1-1/p}[\log(1-x^q)]^2}{qx^{q-1}}-\int\limits_0^x\frac{\log(1-t^q)(\log(1-t^q)-2p)}{(1-t^q)^{1/p}}dt<0.
$$
We have $G(0)=0$. Taking the derivative in $x$ and rearranging, we obtain
$$
G'(x)=\frac{p\log(1-x^q)}{qx^{q}(1-x^q)^{1/p}}\big\{\log(1-x^q)-(q+x^q)\log(1-x^q)\big\}
$$
We want to prove that $G'(x)<0$. The expression in front of the parenthesis is negative. Therefore, we need to show that the expression in parenthesis is positive. But
$$
\log(1-x^q)-(q+x^q)\log(1-x^q)=\log(1-x^q)(1-q-x^q)>0,
$$
and we have proved that $G'(x)<0$. Therefore, $p\to\sin_{p,q}(y)$ is concave.

(ii) In order to prove concavity of $q\to\sin_{p,q}(y)$ we need to show (\ref{eq:qsinpq}). If we drop the first term in the above
inequality, the expression on the left becomes smaller. Hence, if
we can prove that it is still positive, we are finished. This amounts to
showing that
$$
\frac{f}{g}:=\left(\int_0^x\frac{t^q(t^q+p)[\log(1/t)]^2}{(1-t^q)^{2+1/p}}dt\right)/\left(\int_0^x\frac{t^q\log(1/t)dt}{(1-t^q)^{1+1/p}}\right)>\frac{2x^q\log(1/x)}{(1-x^q)}.
$$
We have
$$
\frac{f'}{g'}=\frac{(x^q+p)\log(1/x)}{1-x^q}.
$$
It is easy to check by differentiation  that the function on the right decreases on $(0,1)$.  Clearly, $f(0)=g(0)=0$ and $gg'>0$,
so that we are in the position to apply Lemma~\ref{lm:LHospital}, yielding
$$
\left(\int_0^x\frac{t^q(t^q+p)[\log(t)]^2}{(1-t^q)^{2+1/p}}dt\right)/\left(\int_0^x\frac{t^q\log(t)dt}{(1-t^q)^{1+1/p}}\right)>\frac{(x^q+p)\log(1/x)}{1-x^q}>\frac{2x^q\log(1/x)}{(1-x^q)}.~\square
$$

\medskip

\begin{Theorem}\label{th:sinhlogconvex}
The function $p\to\sinh_{p,q}(y)$ is log-convex on $(1,\infty)$ for any fixed $y\in[0,\pi_q](=\cap_{p>1}[0,\hat{\pi}_{p,q}/2])$ and any fixed $q\in(1,\infty)$. 
\end{Theorem}

\textbf{Proof.}  We need to show that for any $x\in(0,\infty)$ and any $p,q\in(1,\infty)$,
$$
\frac{(q-p)x^q-p}{px(1+x^q)^{1-1/p}}\left(\int_0^x\lambda_1(p,t)dt\right)^2+2\log(1+x^q)\int_0^x\lambda_1(p,t)dt
+\int_0^x\lambda_{2}(p,t)dt\ge0,
$$
where
\begin{align*}
&\lambda_1(p,t)\!=\!\frac{\log(1+t^q)}{(1+t^q)^{1/p}},
\\
&\lambda_2(p,t)\!=\!\frac{(2p-\log(1+t^q))\log(1+t^q)}{(1+t^q)^{1/p}}=2p\lambda_1(p,t)-\frac{\log^2(1+t^q)}{(1+t^q)^{1/p}}.
\end{align*}

Substituting the second expression for $\lambda_2(p,t)$, writing $\mu(p,x)=\int_0^x\lambda_1(p,t)dt$ and rearranging the first terms,
we can rewrite the required inequality as \vspace{-5pt}
\begin{multline*}
\frac{qx^{q-1}}{p(1+x^q)^{1-1/p}}\mu(p,x)^2+2(\log(1+x^q)+p)\mu(p,x)
\\
\ge\frac{(x^q+1)^{1/p}}{x}\mu(p,x)^2+\int_0^x\frac{\log^2(1+t^q)}{(1+t^q)^{1/p}}dt,
\end{multline*}
or, dividing by $\mu(p,x)$,
\begin{multline}\label{eq:divided}
\frac{qx^{q-1}}{p(1+x^q)^{1-1/p}}\mu(p,x)+2(\log(1+x^q)+p)
\\
\ge\frac{(x^q+1)^{1/p}}{x}\mu(p,x)+\frac{1}{\mu(p,x)}\int_0^x\frac{\log^2(1+t^q)}{(1+t^q)^{1/p}}dt.
\end{multline}

Next, we will show that
\begin{equation}\label{eq:srlog}
\frac{1}{\mu(p,x)}\int_0^x\frac{\log^2(1+t^q)}{(1+t^q)^{1/p}}dt
=\left.\underbrace{\int_0^x\frac{\log^2(1+t^q)}{(1+t^q)^{1/p}}dt}_{s(x)}\middle/
\underbrace{\int_0^x\frac{\log(1+t^q)}{(1+t^q)^{1/p}}dt}_{r(x)}<\log(1+x^q)\right..
\end{equation}
Indeed,  $r(0)=s(0)=0$, $s(x)s'(x)>0$ on $(0,\infty)$ and 
$$
\frac{r'(x)}{s'(x)}=\frac{1}{\log(1+x^q)} 
$$
is decreasing.  According to Lemma~\ref{lm:LHospital},
$$
\frac{r(x)}{s(x)}>\frac{r'(x)}{s'(x)}=\frac{1}{\log(1+x^q)},
$$
which is precisely \eqref{eq:srlog}. Hence, if we substitute the rightmost term in \eqref{eq:divided} by $\log(1+x^q)$, we increase the right-hand side, so that the new inequality 
$$
\frac{qx^{q-1}}{p(1+x^q)^{1-1/p}}\mu(p,x)+2(\log(1+x^q)+p)
\ge\frac{(x^q+1)^{1/p}}{x}\mu(p,x)+\log(1+x^q)
$$
is stronger than \eqref{eq:divided}.  Our next goal is to establish this inequality.  Rearranging slightly, we obtain an equivalent form 
\begin{multline}\label{eq:form1}
f_1(x)=\frac{qx^{q}}{p(1+x^q)}\int_0^x\frac{\log(1+t^q)}{(1+t^q)^{1/p}}dt+\frac{x\log(1+x^q)+2px}{(x^q+1)^{1/p}}
\\
\ge\int_0^x\frac{\log(1+t^q)}{(1+t^q)^{1/p}}dt=f_2(x).
\end{multline}
Note first that if 
$$
\frac{qx^{q}}{p(1+x^q)}\ge1,
$$
then inequality \eqref{eq:form1} is trivially true. Hence, we assume that 
\begin{equation}\label{eq:xpqcondition}
0<\frac{qx^{q}}{p(1+x^q)}<1~\Leftrightarrow~qx^{q}-p(1+x^q)<0.
\end{equation}
Clearly, $f_1(0)=f_2(0)=0$.  Differentiating \eqref{eq:form1} with respect to $x$ after substantial simplification, we obtain
$$
f_1'(x)-f_2'(x)=\frac{q^2x^{q-1}}{p(1+x^q)^2}\int_0^x\frac{\log(1+t^q)}{(1+t^q)^{1/p}}dt+\frac{2p(1+x^q)-qx^q}{(1+x^q)^{1+1/p}},
$$
so that $f_1'(x)-f_2'(x)>0$ is equivalent to 
$$
\frac{q^2x^{q-1}}{p(1+x^q)^{1-1/p}}\int_0^x\frac{\log(1+t^q)}{(1+t^q)^{1/p}}dt>qx^q-p(1+x^q)-p(1+x^q)
$$
which is trivially true in view of \eqref{eq:xpqcondition}.  This proves \eqref{eq:form1} and we are done. $\hfill\square$

\section{Evaluation of Some Integrals}\label{sec3}
In this section, the standard notation ${}_pF_q$ is used to denote the generalized hypergeometric function; see \cite{Andrews(1999)} (2.1.2).
In the following theorem, we apply Feynman's trick to evaluate the integrals of the generalized inverse trigonometric functions. In the subsequent corollary, we rewrite these integrals as integrals of certain combinations of direct functions. Note that somewhat similar but simpler integrals (having one free parameter less than ours) have been evaluated in  \cite{Kobayashi(2019)} (Section~3).
\begin{Theorem}\label{th:arcint}
    Suppose $\alpha<2$, $\beta<1$, $0\le{s}\le1$. Then, the following integral evaluation holds:
		\begin{equation}\label{eq:arcsin_int}
			\int_{0}^{1}\frac{\arcsin_{p,q}(sx)}{x^{\alpha}(1-x^q)^{\beta}}dx
			=\frac{s}{q}B\Big(\frac{2-\alpha}{q},1-\beta\Big){_{3}F_{2}}\left(\left.\!\!\begin{array}{c}1/p,(2-\alpha)/q,1/q\\(2-\alpha)/q+1-\beta,1/q+1\end{array}\right|s^q\!\right).
		\end{equation}
If $\alpha,\beta<1$, $0\le{s}\le1$,  then the following integral evaluation holds:
		\begin{multline}\label{eq:arccos_int}
			\int_{0}^{1}\frac{\arccos_{p,q}(sx)}{x^{\alpha}(1-x^p)^{\beta}}dx=\frac{1}{pq}B\Big(\frac{p-1}{p},\frac{1}{q}\Big)B\Big(\frac{1-\alpha}{p},1-\beta\Big)
		\\
		-\frac{s^{p-1}}{q(p-1)}B\Big(1-\frac{\alpha}{p},1-\beta\Big){_{3}F_{2}}\left(\left.\!\!\begin{array}{c}1-1/q,1-\alpha/p,(p-1)/p\\2-\alpha/p-\beta,(2p-1)/p\end{array}\right|s^p\!\right).
		\end{multline}
If $q\alpha>-1$, $q(\alpha-\beta)-\min(1,q/p)<-1$,  $0\le{s}\le1$, then the following integral evaluation holds:
		\begin{multline}\label{eq:arcsinh_int}
			\int_{0}^{\infty}\frac{\arcsinh_{p,q}(sx)}{x^{1-\alpha{q}}(1+x^q)^{\beta}}dx=
			\frac{s}{q}B\Big(\alpha+\frac{1}{q},\beta-\alpha-\frac{1}{q}\Big){_{3}F_{2}}\left(\left.\!\!\begin{array}{c}1/p,\alpha+1/q,1/q\\\alpha-\beta+1/q+1,1/q+1\end{array}\right|s^q\!\right)
			\\
			+\frac{s^{q(\beta-\alpha)}}{q^2(\beta-\alpha)}B\Big(\alpha-\beta+\frac{1}{q},\beta-\alpha+\frac{1}{p}-\frac{1}{q}\Big){_{3}F_{2}}\left(\left.\!\!\begin{array}{c}\beta,\beta-\alpha+1/p-1/q,\beta-\alpha\\\beta-\alpha-1/q+1,\beta-\alpha+1\end{array}\right|s^q\!\right).
		\end{multline}	
If $p<q$, $s\ge1$, $\beta<1$ and $(\alpha+\beta)p>1$, then the following integral evaluation holds:
\begin{multline}\label{eq:arccosh_int}
			\int_{1}^{\infty}\frac{\arccosh_{p,q}(sx)}{x^{\alpha{p}}(x^p-1)^{\beta}}dx=\frac{1}{pq}B\Big(\frac{1}{p}-\frac{1}{q},\frac{1}{q}\Big)B\Big(1-\beta,\beta+\alpha-\frac{1}{p}\Big)
   \\
   +\frac{s^{p/q-1}}{p-q}B\Big(\beta+\alpha-\frac{1}{q},1-\beta\Big){_{3}F_{2}}\left(\left.\!\!\begin{array}{c}1-1/q,\beta+\alpha-1/q,1/p-1/q\\\alpha+1-1/q,1/p-1/q+1\end{array}\right|s^{-p}\!\right).
		\end{multline}	
\end{Theorem}
\textbf{Proof.} The conditions for convergence of the integral in \eqref{eq:arcsin_int} follow from the asymptotic approximations 
 $$
     \arcsin_{p,q}(x)\sim x~~\text{as}~x\to0,~~~~
     \arcsin_{p,q}(x) \sim \pi_{p,q}/2~~\text{as}~x\to1.
 $$ 

Denote
	$$
	A(s)=\int_{0}^{1}\frac{\arcsin_{p,q}(sx)}{x^{\alpha}(1-x^q)^{\beta}}dx.
	$$
Differentiating with respect to $s$ under the integral sign in view of \eqref{eq:arcsinpq} and Euler's integral for the hypergeometric function ${}_2F_1$ \cite{Andrews(1999)} (Theorem~2.2.1), we obtain
	\begin{multline*}
		A'(s)=\int_{0}^{1}\frac{xdx}{x^{\alpha}(1-x^q)^{\beta}(1-s^qx^q)^{1/p}}
		=\frac{1}{q}\int_{0}^{1}\frac{t^{1/q(1-\alpha)}t^{1/q-1}dt}{(1-t)^{\beta}(1-s^qt)^{1/p}}
		\\
		=\frac{B((2-\alpha)/q,1-\beta)}{q}  {_{2}F_{1}}\left(\left.\!\!\begin{array}{c}1/p,1/q(1-\alpha)+1/q\\1/q(1-\alpha)+1/q+1-\beta\end{array}\right|s^q\!\right).
	\end{multline*}
As $A(0)=0$,  term-wise integration and application of $1/(k+\gamma)=(\gamma)_k/[\gamma(\gamma+1)_k]$ yield:
	$$
	A(s)\!=\!\int_{0}^{s}A'(t)dt=\frac{sB((2-\alpha)/q,1-\beta)}{q}\: {}_{3}F_{2}\left(\left.\!\!\begin{array}{c}1/p,1/q(1-\alpha)+1/q,1/q\\1/q(1-\alpha)+1/q+1-\beta,1/q+1\end{array}\right|s^q\!\right).
	$$
	
	The conditions for convergence of the integral in \eqref{eq:arccos_int} follow from the asymptotic approximations (the asymptotics at $x\to1$ is only needed when $s=1$)
 $$
     \arccos_{p,q}(x)\sim  \pi_{p,q}/2~~\text{as}~x\to0,~~~~
     \arccos_{p,q}(x) \sim C(1-x)^{1/q}~~\text{as}~x\to1.
 $$ 
 
Denote
	$$
	G(s)=\int_{0}^{1}\frac{\arccos_{p,q}(sx)}{x^{\alpha}(1-x^p)^{\beta}}dx.
	$$
Differentiating with respect to $s$ under the integral sign in view of \eqref{eq:arccospq} and  Euler's integral for the hypergeometric function ${}_2F_1$ \cite{Andrews(1999)} (Theorem~2.2.1), we obtain
	\begin{multline*}
		G'(s)=-\frac{p}{q}s^{p-2}\int_{0}^{1}\frac{x^{p-1}dx}{x^{\alpha}(1-x^p)^{\beta}(1-s^px^p)^{1-1/q}}
		=-\frac{s^{p-2}}{q}\int_{0}^{1}\frac{t^{-\alpha/p}dt}{(1-t)^{\beta}(1-s^pt)^{1-1/q}}
		\\
		=-s^{p-2}\frac{\Gamma(1-\alpha/p)\Gamma(1-\beta)}{q\Gamma(2-\alpha/p-\beta)}  {_{2}F_{1}}\left(\left.\!\!\begin{array}{c}1-1/q,1-\alpha/p\\2-\alpha/p-\beta\end{array}\right|s^p\!\right).
	\end{multline*}
By changing the variable $t=x^p$, we easily obtain:
	$$
	G(0)=\arccos_{p,q}(0)\int_{0}^{1}x^{-\alpha}(1-x^p)^{-\beta}dx=\frac{1}{pq}B\Big(\frac{p-1}{p},\frac{1}{q}\Big)B\Big(\frac{1-\alpha}{p},1-\beta\Big).
	$$
Hence, we obtain by term-wise integration and an application of $1/(k+\gamma)=(\gamma)_k/[\gamma(\gamma+1)_k]$:
	\begin{multline*}
		G(s)=G(0)+\int_{0}^{s}G'(t)dt=\frac{1}{pq}B\Big(\frac{p-1}{p},\frac{1}{q}\Big)B\Big(\frac{1-\alpha}{p},1-\beta\Big)
		\\
		-\frac{s^{p-1}}{q(p-1)}B(1-\alpha/p,1-\beta){_{3}F_{2}}\left(\left.\!\!\begin{array}{c}1-1/q,1-\alpha/p,(p-1)/p\\2-\alpha/p-\beta,(2p-1)/p\end{array}\right|s^p\!\right).
	\end{multline*}
	
The conditions for convergence of the integral in \eqref{eq:arcsinh_int} follow from the asymptotic approximations 
\begin{align*}
 &\arcsinh_{p,q}(x)\sim  x~~\text{as}~x\to0~~\arcsinh_{p,q}(x) \sim \hat{\pi}_{p,q}/2~~\text{as}~x\to\infty~\text{if}~p<q,\\
     &\arcsinh_{p,q}(x) \sim x^{1-q/p}~~\text{as}~x\to\infty~\text{if}~p>q,
      ~~\arcsinh_{p,q}(x) \sim \log(x)~~\text{as}~x\to\infty~\text{if}~p=q.
\end{align*}

Denote
	$$
	C(s)=\int_{0}^{\infty}\frac{\arcsinh_{p,q}(sx)}{x^{1-\alpha{q}}(1+x^q)^{\beta}}dx.
	$$
Differentiating with respect to $s$ under the integral sign in view of \eqref{eq:arcsinhpq}, we obtain
	\begin{multline*}
		C'(s)=\int_{0}^{\infty}\frac{x^{\alpha{q}}}{(1+x^q)^{\beta}(1+s^qx^q)}dx
		=\frac{1}{q}\int_{0}^{\infty}\frac{t^{\alpha+1/q-1}dt}{(1+t)^{\beta}(1+s^qt)}
		\\
		=\frac{\Gamma(\alpha+1/q)\Gamma(\beta-\alpha+1/p-1/q)}{q\Gamma(\beta+1/p)}  {_{2}F_{1}}\left(\left.\!\!\begin{array}{c}1/p,\alpha+1/q\\\beta+1/p\end{array}\right|1-s^q\!\right),
	\end{multline*}
	where the last equality is an application of the following guise of the Euler integral representation:
	$$
	{_{2}F_{1}}\left(\left.\!\!\begin{array}{c}a,b\\c\end{array}\right|z\!\right)=\frac{\Gamma(c)}{\Gamma(c-b)}\int_{0}^{\infty}\frac{t^{b-1}dt}{(1+t)^{c-a}(1+(1-z)t)^{a}}.
	$$
It is obtained by the change of variable  $t\to{t/(t+1)}$ in \cite{Andrews(1999)} (Theorem~2.2.1). Next, we apply the connection formula \cite{Andrews(1999)} (2.3.13), leading to 
	\begin{multline*}
		C'(s)=
\frac{\Gamma(\alpha+1/q)\Gamma(\beta-\alpha-1/q)}{q\Gamma(\beta)} {_{2}F_{1}}\left(\left.\!\!\begin{array}{c}1/p,\alpha+1/q\\\alpha-\beta+1/q+1\end{array}\right|s^q\!\right)
		\\
		+\frac{\Gamma(\alpha-\beta+1/q)\Gamma(\beta-\alpha+1/p-1/q)}{q\Gamma(1/p)}s^{q(\beta-\alpha)-1} {_{2}F_{1}}\left(\left.\!\!\begin{array}{c}\beta,\beta-\alpha+1/p-1/q\\\beta-\alpha-1/q+1\end{array}\right|s^q\!\right),
	\end{multline*}
As $C(0)=0$, by term-wise integration and an application of $1/(k+\gamma)=(\gamma)_k/[\gamma(\gamma+1)_k]$, we arrive at 
	\begin{multline*}
		C(s)=\int_{0}^{s}C'(t)dt=
		s\frac{\Gamma(\alpha+1/q)\Gamma(\beta-\alpha-1/q)}{q\Gamma(\beta)} {_{3}F_{2}}\left(\left.\!\!\begin{array}{c}1/p,\alpha+1/q,1/q\\\alpha-\beta+1/q+1,1/q+1\end{array}\right|s^q\!\right)
		\\
		+s^{q(\beta-\alpha)}\frac{\Gamma(\alpha-\beta+1/q)\Gamma(\beta-\alpha+1/p-1/q)}{q^2(\beta-\alpha)\Gamma(1/p)} {_{3}F_{2}}\left(\left.\!\!\begin{array}{c}\beta,\beta-\alpha+1/p-1/q,\beta-\alpha\\\beta-\alpha-1/q+1,\beta-\alpha+1\end{array}\right|s^q\!\right).
	\end{multline*}
 
The conditions for convergence of the integral in  \eqref{eq:arccosh_int} follow from the asymptotic approximations under the assumption $p<q$ (the asymptotics at $x\to1$ is only needed when $s=1$):
 $$
     \arccosh_{p,q}(x)\sim  \pi_{p,q}/2~~\text{as}~x\to\infty,~~~~
     \arccosh_{p,q}(x) \sim C(x-1)^{1/q}~~\text{as}~x\to1.
 $$ 
 
Denote
 $$
 D(s)=\int_{1}^{\infty}\frac{\arccosh_{p,q}(sx)}{x^{\alpha{p}}(x^p-1)^{\beta}}dx.
 $$
Differentiating with respect to $s$ under the integral sign in view of \eqref{eq:arccoshpq}, we obtain by substitution $u=x^{-p}$ and Euler's integral for ${}_2F_1$ \cite{Andrews(1999)} (Theorem~2.2.1)
\begin{multline*}
D'(s)=\frac{p}{q}s^{p-2}\int_{1}^{\infty}\frac{x^{p-1}(s^px^p-1)^{1/q-1}}{x^{\alpha{p}}(x^p-1)^{\beta}}dx
=\frac{s^{p/q-2}}{q}\int_{0}^{1}\frac{u^{\beta+\alpha-1/q-1}(1-u)^{-\beta}}{(1-us^{-p})^{1-1/q}}du
\\
=\frac{s^{p/q-2}}{q}B(1-\beta,\beta+\alpha-1/p) {_{2}F_{1}}\left(\left.\!\!\begin{array}{c}1-1/q,\beta+\alpha-1/q\\\alpha+1-1/q\end{array}\right|s^{-p}\!\right).
\end{multline*}
By substitution $u=x^{-p}$ and Euler's beta integral, we see that
$$
D(\infty)=\frac{\hat{\pi}_{p,q}}{2}\int_{1}^{\infty}x^{-\alpha{p}}(x^p-1)^{-\beta}dx=\frac{1}{pq}B(1/p-1/q,1/q)B(1-\beta,\beta+\alpha-1/p),
$$ 
and therefore,
$$
 D(s)=D(\infty)-\int_{s}^{\infty}D'(t)dt.
$$
Performing term-wise integration and applying the relation $1/(k+\gamma)=(\gamma)_k/[\gamma(\gamma+1)_k]$, we arrive at \eqref{eq:arccosh_int}.  $\hfill\square$

\begin{Corollary}
Under the convergence conditions specified in Theorem~\ref{th:arcint}, the following integral evaluations hold:
\begin{multline}\label{eq:sin_int}
	\int\limits_{0}^{x}\frac{t\sin_{p,q}^{-\alpha}(t)\cos_{p,q}(t)dt}{(\sin_{p,q}^q(x)-\sin_{p,q}^q(t))^{\beta}}=\frac{1}{q}B\Big(\frac{2-\alpha}{q},1-\beta\Big)\sin_{p,q}^{2-\alpha-q\beta}(x)
	\\
	\times {_{3}F_{2}}\left(\left.\!\!\begin{array}{c}1/p,(2-\alpha)/q,1/q\\(2-\alpha)/q+1-\beta,1/q+1\end{array}\right|\sin_{p,q}^q(x)\!\right),
\end{multline}
\begin{multline}\label{eq:cos_int}
	\int\limits_{x}^{\pi_{p,q}/2}\frac{t\sin_{p,q}^{q-1}(t)\cos_{p,q}^{2-p-\alpha}(t)}{(\cos_{p,q}^p(x)-\cos_{p,q}^p(t))^{\beta}}dt=\frac{\cos_{p,q}^{1-\alpha-p\beta}(x)}{q^2}B\Big(\frac{p-1}{p},\frac{1}{q}\Big)B\Big(\frac{1-\alpha}{p},1-\beta\Big)
		\\
		-\frac{p\cos_{p,q}^{p(1-beta)-\alpha}(x)}{q^2(p-1)}B\Big(1-\frac{\alpha}{p},1-\beta\Big){_{3}F_{2}}\left(\left.\!\!\begin{array}{c}1-1/q,1-\alpha/p,(p-1)/p\\2-\alpha/p-\beta,(2p-1)/p\end{array}\right|\cos_{p,q}^{p}(x)\!\right),
\end{multline}
\begin{multline}\label{eq:sinh_int}
	\int\limits_{0}^{\hat{\pi}_{p,q}/2}\frac{t\cosh_{p,q}(t)\sinh_{p,q}^{\alpha{q}-1}(t)}{(\sinh_{p,q}^q(x)+\sinh_{p,q}^q(t))^{\beta}}dt
 \\
 =\frac{\sinh_{p,q}^{q(\alpha-\beta)+1}(x)}{q}B\Big(\alpha+\frac{1}{q},\beta-\alpha-\frac{1}{q}\Big) {_{3}F_{2}}\left(\left.\!\!\begin{array}{c}1/p,\alpha+1/q,1/q\\\alpha-\beta+1/q+1,1/q+1\end{array}\right|\sinh_{p,q}^q(x)\!\right)
			\\
			+\frac{B(\alpha-\beta+1/q,\beta-\alpha+1/p-1/q)}{q^2(\beta-\alpha)} {_{3}F_{2}}\left(\left.\!\!\begin{array}{c}\beta,\beta-\alpha+1/p-1/q,\beta-\alpha\\\beta-\alpha-1/q+1,\beta-\alpha+1\end{array}\right|\sinh_{p,q}^q(x)\!\right),
\end{multline}
and
\begin{multline}\label{eq:cosh_int}
	\int\limits_{x}^{\hat{\pi}_{p,q}/2}\frac{t\sinh_{p,q}^{q-1}(t)\cosh_{p,q}^{2-p(\alpha+1)}(t)}{(\cosh_{p,q}^p(t)-\cosh_{p,q}^p(x))^{\beta}}dt=\frac{\cosh_{p,q}^{1-p(\alpha+\beta)}(x)}{q^2}
 B\Big(\frac{1}{p}-\frac{1}{q},\frac{1}{q}\Big)B\Big(1-\beta,\beta+\alpha-\frac{1}{p}\Big)
   \\
  + \frac{\cosh_{p,q}^{p/q-p(\alpha+\beta)}(x)}{q(1-q/p)}B\Big(\beta+\alpha-\frac{1}{q},1-\beta\Big){}_{3}F_{2}\left(\left.\!\!\!\begin{array}{c}1-1/q,\beta+\alpha-1/q,1/p-1/q\\\alpha+1-1/q,1/p-1/q+1\end{array}\!\!\right|\!\cosh_{p,q}^{-p}(x)\!\right).
\end{multline}
\end{Corollary}

{\bf{Remark}}
Note that in \eqref{eq:sin_int} and \eqref{eq:cos_int}
$$
\sin_{p,q}^q(x)-\sin_{p,q}^q(t)=\cos_{p,q}^p(t)-\cos_{p,q}^p(x).
$$

\textbf{Proof.}  To establish \eqref{eq:sin_int}, carry out substitution $x=\sin_{p,q}(t)/s$ (so that $dx=\cos_{p,q}(t)dt/s$) in \eqref{eq:arcsin_int} and write $x=\arcsin_{p,q}(s)$ in the resulting expression.   Similarly, to prove \eqref{eq:cos_int}, carry out substitution $x=\cos_{p,q}(t)/s$ (so that $dx=-q\sin_{p,q}^{q-1}(t)\cos_{p,q}^{2-p}(t)dt/(sp)$) in \eqref{eq:arccos_int} and write $x=\arccos_{p,q}(s)$ in the resulting expression.  The proofs for the remaining integrals are similar. $\hfill\square$

\section{Conclusions}

In this paper, we considered the generalized trigonometric and hyperbolic functions with two parameters. A short survey of their basic properties is presented in the introduction and it includes two presumably new representations for $\arccos_{p,q}$ and $\arccosh_{p,q}$ in terms of the Gauss hypergeometric function.  The main body of the paper is concerned with inequalities for the generalized trigonometric and hyperbolic functions, namely monotonicity and convexity in parameters.  We establish monotonicity for all functions where it takes place and (logarithmic) convexity/concavity for $p\to\tamh_{p',p}$ ($p'=p/(p-1)$), $p\to\sin_{p,q}$, $q\to\sin_{p,q}$ and $p\to\sinh_{p,q}$.  The final section is devoted to integral evaluations. We derive explicit formulas in terms of the generalized hypergeometric functions for four integrals of inverse and four related integrals of direct generalized trigonometric and \mbox{hyperbolic functions.}

Finally, let us mention some potential applications of the results of this paper. Applications in approximation theory have been considered in \cite{BL2015}, where the basis properties of the generalized sine function with two parameters were examined.

Very recently the generalized trigonometric functions were found to play an important role in constructing explicit solutions of a nonlinear generalization of the Schr\"{o}dinger equation. In particular, they enter the general solution of the traveling wave ansatz, as explained in \cite{GPPCT2024}.

Other potential applications include using generalized trigonometric and hyperbolic functions as activation functions for neural networks. In particular, the classical hyperbolic tangent function has been used in this capacity in neural network models; see \cite{DCH(2022)} (Section~3).  Adding two parameters $p$ and $q$ will give an additional degree of flexibility to such models, as they play an entirely different (and highly nonlinear) role from the weights of the model.

\paragraph{Acknowledgements.} The work of the second author was supported by the state assignment of the Institute of Applied Mathematics FEBRAS (No.075-00459-24-00).

\end{document}